\documentclass{amsart}
\usepackage{amsmath,amssymb,amsthm}
\newtheorem{thm}{Theorem}

\newtheorem{cor}{Corollary}
\newtheorem*{lappr}{Approximation Lemma}

\newtheorem{prop}{Proposition}
\theoremstyle{definition}
\newtheorem{ex}{Example}

\theoremstyle{remark}

\newtheorem{defin}{Definition}

\begin{document}
\title{Weighted pluripotential theory on complex K\"ahler manifolds}

\author{Maritza M. Branker}
\address{Department of Mathematics, Niagara University, NY 14109}
\email{mbranker@niagara.edu}
\thanks{}
\author{Ma{\l}gorzata Stawiska}
\address{DFG Research Center Matheon\\ Technische Universit\"{a}t Berlin\\
Berlin\\ Germany}
\email{stawiska@priort.math.tu-berlin.de}
\thanks{Preliminary version}
\dedicatory{} 
\begin{abstract}

We introduce a weighted version of the pluripotential theory on
compact K\"{a}hler manifolds developed by Guedj and Zeriahi. We
give the appropriate definition of a weighted pluricomplex Green
function, its basic properties and consider its behaviour under
holomorphic maps. We also establish a generalization of Siciak's
H-principle.
\end{abstract}
\maketitle

\section*{Introduction}

\noindent Recently there has been significant progress in weighted
pluripotential theory on $\mathbb{C}^N$ which was originally
developed in \cite{Si1},\cite{Si2} and generalized to parabolic manifolds in \cite{Ze}.  Specifically, we refer to
\cite{BL}, \cite{Bl1}, \cite{Bl2}, \cite{Bra}, \cite{MS}.
Concurrently, pluripotential theory on a compact K\"{a}hler
manifold $X$ based on quasiplurisubharmonic functions has been
explored in \cite{GZ1}, \cite{GZ2},\cite{Ko} and \cite{HKH}(see also applications  in \cite{Be1}, \cite{Be2}, \cite{BB}).  The
goal of our article is to develop a framework which would allow
for a unified treatment of both generalizations of the classical
theory and would also allow one to create an analog of the psh-homogeneous pluripotential theory  We will start by showing that a weighted pluripotential
theory on $\mathbb{C}^N$ extends naturally to a pluripotential
theory on $\mathbb{CP}^N$ with a suitably modified weight. In turn
this extends to a homogeneous  pluripotential theory in the universal line
bundle over $\mathbb{CP}^N$, whose charts are biholomorphic to $\mathbb{C}^{N+1}$. We will generalize these results to projective algebraic manifolds. \\
We define a weighted pluricomplex Green function on a compact complex manifold $X$ with a K\"{a}hler form $\omega$.The definition is formulated in terms of a mild
function (see Definition 1). However, many results of our theory hold without requiring that $Q$ be mild.
For a mild function $Q$  and  a Borel set $K \subset X$
 the weighted pluricomplex Green function is
\[
V_{K,\omega,Q}= \sup\{\phi \in PSH(X,\omega): \phi \leq Q \mbox{
on }K\}.
\]
Basic properties of $V_{K,\omega,Q}$ are stated and proved in Section \ref{s:weighted}, followed by the extension of the weighted pluripotential theory in $\mathbb{C}^N$ to a suitable weighted pluripotential theory on $\mathbb{CP}^N$.
We obtain more specific results, in particular a generalized Siciak's H-principle and some approximation results,  in the case when $X$ admits a positive line bundle (which by Kodaira's imbedding theorem is equivalent to $X$ being projective algebraic).

The initial motivation for our work was the similarity between
Theorem 2.12 in \cite{Bra} and Theorem 1 in \cite{St1}; both of
which are generalized versions of Theorem 5.3.1 in \cite{Kl}. We
succeeded in proving the following result (Theorem 5,  Section 2)
which gives the above mentioned theorems as special cases.

\textit{Theorem: Let $(X,\omega)$ be a compact complex K\"{a}hler manifold and $f:X \to X$ a holomorphic surjection. Assume that there exist $\alpha$ and $\beta$, $1<\alpha\leq \beta$,
such that
\newline
$\alpha f_*(PSH(X,\omega))\subset PSH(X,\omega)$ and
$f^*(PSH(X,\omega))\subset \beta\cdot PSH(X,\omega).$ Then for
every Borel set $K\subset X$ and every mild function $Q$ on $X$,
\[
\alpha V_{f^{-1}(K),\omega,f^*Q/\alpha}(x) \leq
V_{K,\omega,Q}\circ f(x) \leq \beta
V_{f^{-1}(K),\omega,f^*Q/\beta}.
\]
}

\section{Weighted pluricomplex Green functions}
 \label{s:weighted}

Throughout the paper we assume that $X$ is a connected compact
complex K\"{a}hler manifold. Therefore we have on $X$ (cf.
\cite{GF}, VI.3) the fundamental form $\omega$ of a hermitian
metric $\Gamma$ on $X$ with $\omega =
i\sum_{j,k}\gamma_{jk}dz_j\wedge d\overline{z}_k$, satisfying
$d\omega=0$. It follows that in each coordinate neighborhood in
$X$ we can define a $\mathcal{C}^{\infty}$ real-valued function
$\phi$ such that $i\partial\overline{\partial}\phi = (1/2)dd^c
\phi = \omega$. The functions $\phi$ are called local potentials
of the K\"{a}hler metric $\Gamma$. Existence of smooth local
potentials is in fact equivalent to $\Gamma$ being K\"{a}hler: if
the fundamental form $\omega$ of $\Gamma$ satisfies $\omega =
i\partial\overline{\partial}\phi$, then $d\omega =0$. For example,
the Fubini- Study metric on $\mathbb{CP}^N$ is K\"{a}hler, since
it has  local potentials given by $\phi_j= \log (1 + \sum_{k \ne
j}|z_{j,k}|^2)$ in the coordinate neighborhoods $U_j = \{Z_j \ne
0\}$ with $j = 0,1,...,N$. Here $[Z_0:...:Z_N]$ are homogeneous
coordinates in $\mathbb{CP}^N$ and $z_{j,k}:=Z_k/Z_j$ in $U_j$.
The set $U_0$ is identified with $\mathbb{C}^N$ and $z_{0,k}=:z_k,
\quad k=1,...,N$, are affine coordinates. We have $\phi_0 = \log(1
+ \|z\|^2)$ for $z \in \mathbb{C}^N$.
 Let $\omega$ be a closed real $(1,1)$ current on $X$
with continuous local potentials. From  \cite{GZ1}, the class of
$\omega$-plurisubharmonic functions  is defined as
\[
PSH(X,\omega) = \{v \in L^1(X,\mathbb{R}\cup\{-\infty\}): dd^c v
\geq -\omega \mbox{ and }v \mbox{ is upper semicontinuous}\}.
\]
 The $\omega$-pluricomplex Green function of a
Borel set $K \subset X$ is defined as
\[
V_{K,\omega}(x) = \sup\{v(x): v \in PSH(X,\omega): v \mid_K \leq
0\}
\]
Consider the class of $PSH(X,\omega)$, where $\omega$ is a
K\"{a}hler form on $X$ with local potentials $\phi_j:U_j\mapsto
\mathbb{R}$, where $\{U_j\}_{j=0}^m$ is an open cover of $X$ by
coordinate neighborhoods.

\begin{defin}\label{def:mild} Let $Q:X\mapsto \mathbb{R}\cup\{+\infty\}$ be a function
 such that $\exp(-Q+\phi_j)$ is
continuous in $U_j, \quad j=1,...,m$ and $\{Q\neq +\infty\}$ is
not a pluripolar subset of $X$. We will call $Q$ satisfying these
assumptions a mild function.
\end{defin}
Mild functions are necessarily lower semicontinuous.

\begin{defin}
For a mild function $Q$ on $X$ and for a Borel set $K \subset X$
let us define the weighted $\omega$-pluricomplex Green function as
\[
V_{K,\omega,Q}= \sup\{\phi \in PSH(X,\omega): \phi \leq Q \mbox{
on }K\}.
\]
\end{defin}

The following properties are direct consequences of our definition
of $V_{K,\omega,Q}$.

\begin{prop}\label{prop:obvious}
Let $K, K_1,K_2$ be Borel subsets of $X$ and $Q, Q_1,Q_2$ be mild functions.
\begin{enumerate}
\item[i)]If $Q_1 \leq Q_2$ on $K$ then $V_{K,\omega,Q_1} \leq V_{K,
\omega, Q_2}$.
\item[ii)]If $K_1 \subset K_2$ then $V_{K_2,\omega, Q} \leq V_{K_1,
\omega, Q}$.
\item[iii)] Let $Q$ be a mild function that belongs to the class
$PSH(X,\omega)$. Then $V_{X,\omega,Q} = Q$.
\item[(iv)] Let $\omega'$ be
cohomologous to $\omega$, $\omega'=\omega + dd^c\xi$ for $ \xi \in
L^1(X)$. If $\xi$ is mild and continuous, then $V_{K,\omega',Q}=
V_{K,\omega,Q-\xi}+\xi$.
\end{enumerate}
\end{prop}
We continue to establish basic properties of the weighted pluricomplex Green function in Propositions 2 and 3.
\begin{prop}\label{prop:nonpolar}
Let $K$ be a Borel set in $X$ and $Q$ a mild function on $X$. If $K$ is
not $PSH(X,\omega)$-polar then $V^{*}_{K,\omega,Q} \in
PSH(X,\omega)$.
\end{prop}
\begin{proof}
By Choquet's lemma there exists an increasing sequence of
functions
\newline
 $\phi_j \in PSH(X, \omega)$ such that $\phi_j \leq Q$ on
$K$ and $$V^{*}_{K,\omega, Q} = (\lim_{j\rightarrow \infty}
\phi_j)^*.$$ It follows from Proposition 2.6(2) in \cite{GZ1} that
$V_{K,\omega,Q}^* \in PSH(X,\omega)$.
\end{proof}

\begin{prop}
Let $E$ be a Borel subset of $X$ and $P$ a $PSH(X,\omega)$ polar
set. Then we have $$V^*_{E \cup P, \omega, Q} =
V^*_{E,\omega,Q}.$$
\end{prop}
\begin{proof}
Recall that a set $P$ is said to be $PSH(X,\omega)$-polar if it is
included in the $- \infty$-locus of some function $\psi \in PSH(X,
\omega)$ which is not identically $-\infty$ on $X$. By Prop
\ref{prop:obvious}(ii) we have $V^{*}_{E \cup P, \omega, Q}
\leq V^{*}_{E,\omega,Q}$. Consequently, we need to establish
\newline
$V^{*}_{E,\omega,Q} \leq V^{*}_{E \cup P, \omega, Q}$. Suppose $u \in
PSH(X, \omega)$ with $u \leq Q$ on $E$ and let $v \in PSH(X, \omega)$ such
that $ P \subset \{ v= - \infty\}$. We may assume $v \leq Q$ on
$E$. Then for each $\epsilon > 0$,
$$(1-\epsilon) u + \epsilon v \leq V_{E \cup P, \omega, Q} \leq V^{*}_{E \cup P, \omega, Q}.$$
Therefore $u \leq V^{*}_{E \cup P, \omega, Q}$ on $X$ and by
taking the supremum,  $V^{*}_{E,\omega,Q} \leq V^{*}_{E \cup P,
\omega, Q}$.

\end{proof}

Let us now show how weighted pluripotential theory on
$\mathbb{C}^N$ can be extended to a suitable weighted
pluripotential theory on $\mathbb{CP}^N$. Recall that in
the weighted theory on $\mathbb{C}^N$ one begins with
an admissible weight function on a closed set $K \subset
\mathbb{C}^N$. An admissible weight $w$ is a nonnegative upper
semicontinuous function $w$ on $\mathbb{C}^N$ with $\{z \in K:
w(z)
> 0\}$ non-pluripolar and satisfying the boundedness condition $\lim_{|z| \to \infty}|z|w(z)=0$ if $K$ is an unbounded set
(cf. \cite{BL}, \cite{Bl1}, \cite{ST}).
 The weighted pluricomplex Green function of $K$
is defined as
\[
V_{K,Q}= \sup\{u \in \mathcal{L}, u \leq Q \mbox{ on }K\}.
\]
where $Q= -\log w$.\\
In the homogeneous coordinates $[Z_0:...:Z_N]$ on $\mathbb{CP}^N$
(with the usual identification $\mathbb{C}^N\simeq \{Z_0 \ne 0\}$
and affine coordinates $z_j = Z_j/Z_0, \quad j=1,...,N$) let
$\tilde{w}(Z_0:...:Z_N)=w(z_1,...,z_N)/|Z_0|$ in $\{Z_0\ne 0\}$,
where $w$ is nonnegative and upper semicontinuous with $\{w>0\}$
non-pluripolar, but not necessarily satisfying the boundedness
condition. The expression $W(Z)=: \|Z\|\tilde{w}(Z)$ defines a
homogeneous function of order 0 in $\mathbb{C}^{N+1}\setminus
\{Z_0 =0\}$. We have $W(Z)=\varphi_0(z)+\log w(z)$ for $Z_0 \neq 0$, where $\varphi_(z)=(1/2)\log(1+|z|^2)$.
We take $$
\sqrt{|Z_1|^2+...+|Z_N|^2}\tilde{w}(0:Z_1:...:Z_N)=\limsup_{0 \ne
Y_0 \to 0,Y_j \to Z_j}\|Y||\tilde{w}(Y), Y = (Y_0,...,Y_N)$$ to
obtain an upper semicontinuous function (still denoted by $W$) globally on
$\mathbb{CP}^N$, with all values greater or equal to 0.
The
boundedness condition is equivalent  to the property that this
global function is identically zero on the hyperplane $\{Z_0=0\}$
. This is because $\lim_{|z|\to \infty}|z|w(z)=\lim_{|z|\to
\infty}\sqrt{1+|z|^2}w(z)$.  We will assume a weaker condition, namely that
$W$ is bounded in $\mathbb{CP}^N$.

The following example demonstrates
that the boundedness condition is too restrictive when
constructing a weighted pluripotential theory on complex
manifolds.
\begin{ex}\label{ex:compact} Let $\omega_{FS}$ be the
Fubini-Study K\"{a}hler form on $X=\mathbb{CP}^N$ with local
potentials $\phi_j$ as above and let $K$ be a  subset of
$\mathbb{C}^N \subset \mathbb{CP}^N$. For $Z \in \mathbb{CP}^N$
define $Q_j(Z)=\phi_j(Z),\quad
 j=0,...,N$ , so that
$Q_0(z)= (1/2)\log (\sqrt{1+\|z\|^2})$ for $ z \in \mathbb{C}^N$.
The natural 1-to-1 correspondence between $PSH(X,\omega_{FS})$ and
the class $\mathcal{L}(\mathbb{C}^N)$ of plurisubharmonic
functions with logarithmic growth at infinity, presented
explicitly in Example 1.2 in
\cite{GZ1}, gives the following:\\
\begin{align*}
V_{K,Q_0}(x) &= \sup \{u(x):u \in \mathcal{L}(\mathbb{C}^N),  u(z)
\leq \log \sqrt{1 + \|z\|^2} \quad \forall z \in K \} \\&= \sup \{u(x):
u \in \mathcal{L}(\mathbb{C}^N),u(z)- (1/2)\log(1+|z|^2)\leq
0 \quad \forall z \in K \}\\&= \sup \{v(x)+ (1/2)\log(1+|x|^2),v \in
PSH(\mathbb{CP}^N,\omega_{FS}): \quad v\mid_K \leq 0 \}
\\&=  V_{K,\omega_{FS}}(x)+ (1/2)\log(1+|x|^2)
\end{align*}
for every $x\in \mathbb{C}^N$. Assume  now that $K$ is not
$PSH(\mathbb{CP}^N,\omega_{FS})$-polar. Then $V_{K,\omega_{FS}}^*
\in PSH(\mathbb{CP}^N,\omega_{FS})$ and $V_{K,Q_0}^* \in
\mathcal{L}(\mathbb{C}^N)$. For a point $Z$ on the hyperplane at
infinity $\{Z_0=0\}$ we get $$V_{K,\omega_{FS}}^*(Z) = \limsup_{x
\to Z, x \in \mathbb{C}^N}(V_{K,Q_0}^*(x)-(1/2)\log(1+|x|^2)).$$
 \end{ex}
  Note that the function $w(z) = \exp(-Q_0(z))$ in our example does not
satisfy the boundedness condition in $\mathbb{C}^N$. Indeed, the
function
 $\|Z\|\tilde{w(z)}=\exp(-Q_j(Z)+\phi_j(Z))$ for $ Z \in U_j, \quad j=0,...,N$ is a
 constant function 1 on $\mathbb{CP}^N$ (which of course is continuous, but never 0).
We draw the reader's attention to the paper \cite{Bl2}, in which a
relation between weighted theory in $\mathbb{C}^N$ and standard
pluripotential theory in $\mathbb{C}^{N+1}$ is outlined.  Examples
considered in the Section 5 of that paper deal with a weight
function $w$ which is given as the Hartogs radius of a domain with
balanced fibers in $\mathbb{C}^{N+1}$ (for the definition and
basic properties, see \cite{Sh}). Such a function is upper
semicontinuous, but as shown in \cite{Bl2}, does not have to
satisfy the boundedness condition on $\mathbb{C}^N$. Furthermore,
the results of \cite{Si2} as well as \cite{MS} were obtained
without assuming the boundedness condition.  It thus seems
reasonable to weaken this condition when working on complex
manifolds. In \cite{Gu} a notion of a 'convex' hull with respect
to a closed real $(1,1)$-current $T$ is considered where the
functions $f$ defining the hull satisfy the condition that
$\exp(f+\phi)$ are continuous, with $\phi$ continuous local
potentials for $T$. We adopted an analogous  condition as a
part of our definition of a mild function.\\
\newline
The method demonstrated in Example \ref{ex:compact}
can also be used to prove the following:
\begin{prop}\label{prop:extend} Let $K \subset \mathbb{C}^N\cong \{Z_0\ne 0\}$.
For a mild function $Q$ on $\mathbb{CP}^N$ with respect to $\omega
= \omega_{FS}$ define
$$q(z_1,...,z_N) = q(Z_1/Z_0,...,Z_N/Z_0)=Q(Z)-\log(\|Z\|/|Z_0|), \quad
Z_0 \ne 0.$$ Conversely, for a lower semicontinuous $q$ on
$\mathbb{C}^N$, consider
$$Q(Z)=q(Z_1/Z_0,...,Z_N/Z_0)+\log\|Z\|+\log|Z_0|,$$
together with its lower semicontinuous regularization as $Z_0 \to
0$. Then $V_{K,q}(x) =.V_{K,\omega,Q}(x)=(1/2)\log(1+\|x\|^2), \quad x \in \mathbb{C}^N$.
\end{prop}

Consider now a holomorphic line bundle $L$ over a compact
K\"{a}hler manifold $X$. Recall that a (singular) metric on $L$
can be given (cf. \cite{De}, \cite{DPS}) by a collection of
real-valued functions $h=\{h_j\}$ on $X$, defined in a
trivializing cover $\{U_j\}$, such that $h_j = h_i+\log|g_{ij}|$,
where $g_{ij}$ are transition functions for $L$. The metric is
called positive if all $h_j$ are plurisubharmonic. (The notion of
positivity is used here in the weak sense.) In particular, a
smooth metric $\{\phi_j\}$ such
that $\omega = dd^c \phi_j$ is a K\"{a}hler form will be positive.\\

If  $L$ is a positive line bundle and $\omega =[c_1(L)]$, there is
a 1-to-1 correspondence between the family of all positive metrics
on $L$ and the class $PSH(X,\omega)$. In the case
of $X=\mathbb{CP}^N$ with the Fubini-Study form $\omega$, this correspondence is equivalent to
the $H$-principle due to Siciak (\cite{Si3}).

\begin{prop}\label{prop:homogeneous}(cf.\cite{G}, property (iv) pg 456): Let $h$ be a logarithmically
homogeneous plurisubharmonic nonnegative function on
$\mathbb{C}^{N+1}$. Then $h$ defines a positive metric on
$\mathbb{CP}^N$. Conversely, each positive metric on
$\mathbb{CP}^N$ defines a logarithmically homogeneous psh function
on $\mathbb{C}^{N+1}$.
\end{prop}
\begin{proof}
By logarithmic homogeneity we have,
$$v(Z_0/Z_k,...,1,...,Z_N/Z_k)=v(Z)-\log|Z_k| \text{ in } \{Z_k
\ne 0\}$$ Hence $v_k = v_i+\log|Z_k/Z_i|$ in $U_i\cap U_k$ and all
$v_i$ are plurisubharmonic. To prove the converse, take
$h_0=h\mid_{U_0}$. The function $v(Z)=h_0(Z)+\log|Z_0|$ in $U_0$,
and $v(0,Z_1,...,Z_N)=\limsup_{(Y_0\to 0, Y_j \to Z_j)}
v(Y_0,Y_1,...,Y_N)$ is plurisubharmonic. Since it also satisfies
$v(\lambda Z)=v(Z)+\log|\lambda|$ for $\lambda \in \mathbb{C}$ our
proof is complete.
\end{proof}
By Example 1.2 in \cite{GZ1}, the class
$\mathcal{L}(\mathbb{C}^N)$ corresponds in a 1-to-1 manner with
the class of $PSH(\mathbb{CP}^N,\omega)$ functions, which in turn
correspond in a 1-to-1 manner with positive metrics on the
(positive) hyperplane bundle over $\mathbb{CP}^N$. Thus
Proposition 5 establishes a 1-to-1 correspondence between
logarithmically homogeneous functions $\tilde{v}$ on
$\mathbb{C}^{N+1}$ and functions $v$ in the class
$\mathcal{L}(\mathbb{C}^N)$ so that $\tilde{v}(1,z)=v(z)$ for $z
\in  \mathbb{C}^N$,that is, the $H$-principle.

If $L$ is a positive line bundle over $X$, then its dual $L'$ is
negative (\cite{GF}, Prop. VI.6.1 and VI.6.2). Hence there exists
a system of trivializations $\theta_i:L'\mid_{U_i}\mapsto U_i
\times \mathbb{C}$ with transition functions $G_{ik} =
g_{ik}^{-1}=g_{ki}$ and a smooth metric $\{h_i\}$ on $L$ such that
the smooth function $\chi_h: L' \mapsto \mathbb{R}$, defined as
$\chi_h\circ \theta_i^{-1}(x,t)= H_i(x)\cdot|t|^2$, is strictly
plurisubharmonic outside the zero section of $L'$, where
$H_i(x)=\exp2h_i(x),\quad x \in U_i$. As a simple example of a
negative line bundle we can take the universal line bundle over
$\mathbb{CP}^N$, $\mathcal{O}(-1) := \{([Z],\xi)\in \mathbb{CP}^N
\times \mathbb{C}^N: \xi \in \mathbb{C}\cdot Z, Z \in
\mathbb{C}^{N+1}\setminus \{0\}, [Z] = \mathbb{C}^*\cdot Z\}$.
That is, the fiber of $\mathcal{O}(-1)$ over a point $[Z] \in
\mathbb{CP}^N$ is the complex line in $\mathbb{C}^{N+1}$ generated
by $(Z_0,...,Z_N)$. The function $\chi\circ\theta_i^{-1}(Z,t)=
|t|^2|Z_i|^{-2}\|Z\|^2$ for  $ Z_i \ne 0$, associated with the
Fubini-Study metric on the dual line bundle $\mathcal{O}(1)$ over
$\mathbb{CP}^N$, is plurisubharmonic.

Next we establish a  generalization of Siciak's $H$-principle.

\begin{thm}\label{thm:H-principle}(cf. \cite{GF}, Prop. VI.6.1):
Let $L$ be a positive line bundle over a compact K\"{a}hler
manifold $X$ and let $d>0$. Let $\mathcal{H}_d^+$ denote the
family of all functions $\chi \in PSH(L')$ which are nonnegative,
not identically 0 and absolutely homogeneous of order $d$ in each
fiber. Then there is a one-to-one correspondence between
$\mathcal{H}_d^+$ and the class of positive metrics on $L$.
\end{thm}
\begin{proof} Consider a system of trivializations $\theta_i:
L'\mid_{U_i}\mapsto U_i\times \mathbb{C}$ with transition
functions $G_{ik}=g_{ki}=1/g_{ik}$. Let $\chi \in
\mathcal{H}_d^+$. For  $x \in U_i,t\ne 0$ define
$$H_i(x):=\chi\circ\theta_i^{-1}(x,t)/|t|^d.$$
Note that this expression does not depend on $t$. We have
$\chi\circ\theta_i^{-1}(x,t)=\chi\circ\theta_k^{-1}(x,G_{ki}(x)t)$,
hence by absolute homogeneity of order d,
$H_k(x)=|G_{ki}(x)|^dH_i(x)$ in $U_i \cap U_k$. Taking $h_i =
(1/d)\log H_i$ in $U_i$ we get a collection of plurisubharmonic
functions satisfying $h_k = \log|g_{ik}| + h_i$, i.e., a positive
metric on $L$. Conversely, let $\{h_i\}$ be a metric on $L$. The
function $\chi$ on $L'$ defined as
$\chi\circ\theta_i^{-1}(x,t)=\exp (dh_i(x))\cdot|t|^d$ is
plurisubharmonic if and only if $h_i$ are, so for a positive
metric the associated function $\chi$ is in $\mathcal{H}_d^+$.
\end{proof}
Unless otherwise indicated, we will work with
$\mathcal{H}^+:=\mathcal{H}_1^+$. Note that if we take $L'$ in
Theorem \ref{thm:H-principle} to be the universal line bundle
$\mathcal{U}$ over $\mathbb{CP}^N$, then the trivialization
$\theta_i:\pi^{-1}(U_i)\mapsto U_i\times \mathbb{C}$ is given as
$\theta_i(t(Z))=([Z_0:...:Z_N],tZ_i)$. Hence for a function $\chi
\in \mathcal{H}^+$ we have $\chi\circ
\theta_i^{-1}([Z_0:...:Z_N],t)=
h_i([Z_0:...:Z_N])+\log|Z_i|+\log|t|$ for $ \quad Z_i \ne 0$,
where $h_i$ define a metric on $\mathbb{CP}^N$. By Proposition
\ref{prop:homogeneous}, over the chart $Z_0\ne 0$ we get
$\chi(tZ)= v(Z_1/Z_0,...,Z_N/Z_0)+\log|t|$ for $ t \ne 0$ with $v$
plurisubharmonic. That is, $\chi$ defines a logarithmically
homogeneous psh function on
$\mathbb{C}^{N+1}$. \\

For a positive holomorphic line bundle $L$ over a compact
K\"{a}hler manifold $X$ there is a precise relation between
the weighted pluricomplex Green function with respect to
$\omega=[c_1(L)]$ of a Borel set K in $X$ and a
$\mathcal{H}^+$-envelope of some associated set $\tilde{K}$ in the
dual bundle $L'$. It  generalizes the formulas obtained by  Bloom in
(\cite{Bl2}).

 For the weight $Q$ on $X$ consider the collection
$q_i=Q-\phi_i$, where $\omega = dd^c\phi_i$ in $U_i$ and $U_i$
form a trivializing cover for $L$. For $K \subset X$ define
$\tilde{K}\subset L'$ by taking $$\tilde{K}\cap
\pi^{-1}(U_i)=\{\theta_i^{-1}(x,t): x \in U_i \cap K, |t| =
\exp(-q_i(x))\}.$$ This set is well defined, since
$\theta_k^{-1}(x,t)=\theta_i^{-1}(x,G_{ki}(x)t)$. Hence if $x \in
U_i \cap U_k \cap K$, then $|G_{ki}(x)t|=\exp(-q_i(x))$ if and
only if $|t|=\exp(-q_k(x))$. Consider $$H_{\tilde{K}}=\sup\{u\in
PSH(L'):\exp u \in \mathcal{H}^+, u\mid_{\tilde{K}}\leq 0\}.$$ The
following theorem gives the relationship between functions
$H_{\tilde{K}}$ and $V_{K,\omega,Q}$.
\begin{thm}\label{thm:preciserel}(cf. \cite{Bl2}, Thm 2.1): For all $i$,
\[
H_{\tilde{K}}\circ\theta_i^{-1}(x,t)=V_{K,\omega,Q}(x)+\log|t|+\phi_i(x)
\]
\end{thm}
\begin{proof} By Theorem \ref{thm:H-principle},
\begin{align*}
H_{\tilde{K}}&=\sup\{u:
u\circ\theta_i^{-1}(x,t)=h_i(x)+\log|t|,\quad u\mid_{U_i\cap
K}\leq 0\} \\ &=\sup\{u:u\circ\theta_i^{-1}(x,t)=h_i(x)+\log|t|,\quad
h_i(x)\leq q_i, \forall i\}
\end{align*}
 where $h_i$ define a positive metric
on $L$. Hence, for such $h_i$,
\begin{align*}
H_{\tilde{K}}\circ\theta_i^{-1}(x,t)&= \sup\{h_i(x):
h_i(x)\mid_{K\cap U_i} \leq q_i\}+\log|t| \\&=
\sup\{v(x)+\phi_i(x):v\in PSH(X,\omega),v\mid_K\leq Q\}+\log|t|\\&
=V_{K,\omega,Q}(x)+\log|t|+\phi_i(x),\quad \forall i.\end{align*}
\end{proof}

Theorem \ref{thm:preciserel} allows us to study the behavior of
the weighted pluricomplex Green functions as we vary the weight.
Namely, we have the following:\\
\begin{prop}\label{prop:increase}(cf. \cite{Bl2}, Cor 2.2):
Let $K \subset X$ be a Borel set. Suppose $\{Q_n\}, Q$ are mild
functions with $Q_n \nearrow Q$. Then $\lim_{n\rightarrow \infty}
V_{K, \omega, Q_n} = V_{K,\omega,Q}.$
\end{prop}
\begin{proof} Consider the
sets $K_n, M_n \subset L'$, where
$$M_n\cap\pi^{-1}(U_i)=\{\theta_i^{-1}(x,t): x \in U_i \cap
K, |t| \leq \exp({-q_i}^{(n)}(x))\},$$

$$K_n\cap\pi^{-1}(U_i)=\{\theta_i^{-1}(x,t): x \in U_i \cap K,
|t| = \exp({-q_i}^{(n)}(x))\}$$ where $ q_i^{(n)}=Q_n-\phi_i, n
\geq 0$. The sequence $M_n$ is decreasing, with
$\bigcap_{n=1}^{\infty}M_n = \tilde{M_0}$.   By maximum principle
(applied in each fiber), $H_{M_n}=H_{K_n}, n \geq 0$ (here we use
the assumption of all $Q_n$ being mild). For a function $u \in
\mathcal{H}^+$ such that $u \leq 0$ on $M_0$ and an arbitrary
$\varepsilon > 0$, there exists an $n_0$ such that for all $n \geq
n_0$ we have $M_n \subset \{u<\varepsilon\}$. The function
$u-\varepsilon$ is in $\mathcal{H}^+$ and for $n \geq n_0$ it
satisfies $u-\varepsilon \leq H_{M_n} \leq \lim_{n \to
\infty}H_{M_n} \leq H_{M_0}$, hence $\lim_{n \to \infty}H_{M_n} =
H_{K_0}$. By Theorem \ref{thm:preciserel}, $\lim_{n \to \infty}
V_{K,\omega,Q_n}=V_{K,\omega,Q_0}$.
\end{proof}

\begin{prop}\label{prop:decrease}(cf. \cite{Bl2}, Cor 2.4) Let $Q_n, n \geq 0$ be mild
functions on $X$ such that $Q_n \searrow Q_0$. Then
$V_{K,\omega,Q_0}= \lim_{n \to \infty}V_{K,\omega,Q_n}$.
\end{prop}
\begin{proof} Since the potentials $\phi_i$ of $\omega$ are
continuous, we have
$H_{\tilde{K}}^*\circ\theta_i^{-1}(x,t)=V_{K,\omega,Q_0}^*+\log|t|+\phi_i$
for all $i$. We can assume that the set $M_1$ (see Proposition
\ref{prop:increase}) is not $\omega$-polar. By Proposition
\ref{prop:nonpolar}, $H_{K_0}^*$ is plurisubharmonic on $L'$. Let
$H=\lim_{n \to \infty}H_{M_n}$. The function $H$ is in
$\mathcal{H}^+$ and satisfies $H \leq 0$ on $K_0\setminus P$,
where $P$ is some pluripolar set. Hence $H \leq H_{K_0}^*$.
\end{proof}
\begin{cor} Proposition \ref{prop:increase},  holds when the convergence
 $Q_n \searrow
Q$, takes place quasi-everywhere on
$X$, that is, outside some $\omega$-polar set.
\end{cor}
\begin{cor} Proposition \ref{prop:decrease} holds when the convergence $Q_n \nearrow Q$ takes place quasi-everywhere on X.
\end{cor}
\section{Approximation and pullbacks by holomorphic maps} \label{s:appr-pull}

In standard pluripotential theory in $\mathbb{C}^N$ and its
weighted generalization there is a function $\Phi_K$ such that
$\log \Phi_K = V_{K,Q}$. The function $\Phi_K$ is given as the
supremum of certain functions with 'regular' growth, that is,
polynomials (when $Q \equiv 0$) or weighted polynomials (see
Theorem 6.2 in\cite{Si1}, Theorem 2.8 in \cite{Bl1}, and
 Th\'{e}or\`{e}me 5.1 in \cite{Ze}). In \cite{GZ1} it is  proven that
$V_{K,\omega}(x)=\sup\{(1/n)\log\|s\|_{n\varphi}(x):n \geq 1, s
\in \Gamma(X,L^n), \sup_K\|s\|_{n\varphi}\leq 1\}$, where $L$
is a positive holomorphic line bundle over a compact manifold
$X$, $\omega = dd^c\varphi_j$  in a trivializing cover $U_j$ is
a (global) K\"{a}hler form and the norm $\|s\|_{n\varphi}$ of a
section $s$ of the tensor power $L^n$ is computed as follows:
$\|s\|_{n\varphi} = |s_j|\exp(-n\varphi_j)$ in $U_j$. All such
theorems are based on the possibility of approximation of
general plurisubharmonic functions by so-called Hartogs
functions, which are obtained by certain operations from
functions of the type $\log|f|$  with $f$ holomorphic (cf.
\cite{Kl}, theorem 5.1.6) Such approximation may be not always
possible, but is possible for example in pseudoconvex domains
in $\mathbb{C}^N$, as shown in \cite{Bre}. Below, we will work
in pseudoconvex neighborhoods of the zero section of $L'$ to
prove the following:

\begin{thm}\label{thm:reggrowth} Let $X,L,\varphi,\omega$ be as above. Let $Q$ be a mild function on
$X$ and let $K$ be a compact subset of $X$. Then\\
$$V_{K,\omega,Q}=\log \Phi_{K,\omega,Q} \text{ where } \Phi_K(x)=\sup_{n\geq
1}(\Phi_n(x))^{1/n}$$ with
$$\Phi_n(x)=\sup\{\|s\|_{n\varphi}(x):n \geq 1, s \in \Gamma(X,L^n),
\sup_K\exp(-nQ)\|s\|_{n\varphi}\leq 1\}.$$
\end{thm}

Unlike \cite{GZ1}, in which the theorem was proved for $Q\equiv 0$, we will not use $L^2$-estimates for the $\bar{\partial}$-operator. Instead, we will apply the following lemma (cf. \cite{Ze}, Lemme
5.2, \cite{Be1}, Lemma 2.1 and 3.2):
\begin{lappr} Let $X,\omega,L$ be as above and let $v \in PSH(X,\omega)\cap \mathcal{C}^\infty$ be such that $dd^cv+\omega$ is strictly positive.
Then for every $\varepsilon > 0$ and every compact $K \subset X$
there exist $N_1,...,N_m$ and $s_1,...,s_m$ such that $s_j \in
\Gamma(X,L^{N_j}),j=1,...,m$ and
$$v(x)-\varepsilon \leq
\sup_{1\leq j \leq m}(1/N_j)\log \|s_j(x)\|_{N_j\varphi} \leq
v(x) \text{ for all } x \in K$$, where the norm of the section
$s_j$ is computed as above.
\end{lappr}

\begin{proof}(of the Approximation Lemma): Let $\varphi_i$ be local
potentials for the K\"{a}hler form $\omega$ and let
$h=\{h_i=v+\varphi_i\}$ be the positive  metric corresponding
to $v$. The inequality in the statement of the lemma is
equivalent to $$h_i-\varepsilon \leq \sup_{1\leq j \leq
m}(1/N_j)\log |s_j(x)| \leq |h_i(x)|, \quad x \in K\cap U_i,
\quad i=1,...,l$$, where $|\cdot|$ is the usual absolute value
of a complex number. Let $r\in (0,1)$ and let $\chi_r$ be the
function in the class $\mathcal{H}^+$ on $L'$ associated with
 the metric $r\cdot h$. For every $r$ the set $ \Omega_r=\{\chi_r<1\}$ is a
strictly pseudoconvex neighborhood of the zero section in $L'$
(cf. \cite{GF}, VI.6.1). Fix a point $x_0 \in K$ and $\zeta_0=
\theta_i^{-1}(x_0,1)$. Then $|t|<\chi_r(\zeta_0)$ if and only
if $(x_0,t) \in\Omega:= \Omega_r$. The function $f(t) =
\sum_{n=1}^{\infty}(\chi_r(\zeta_0))^n t^n, \quad |t| <
1/\chi_r(\zeta_0), \quad f(0) = 0$ is holomorphic on the
analytic set $(\Omega \cap L'_{x_0})\cup X$ and is identically
$0$ on X. Let us consider the Remmert reduction of $\Omega$
(see \cite{G}, Satz 1, or \cite{P}, Theorem 2.7 and preceding
discussion). That is, we have a Stein space $Y$ and a proper
surjective holomorphic map $\Phi:\Omega \mapsto Y$ with the
following properties: (i) $\Phi$ has connected fibers; (ii)
$\Phi_*(\mathcal{O}_\Omega) = \mathcal{O}_Y$;(iii) the
canonical map $\mathcal{O}_Y(Y)\mapsto
\mathcal{O}_\Omega(\Omega)$ is an isomorphism;(iv) if
$\sigma:\Omega \mapsto Z$ is a holomorphic map into a Stein
space $Z$ then there exists a uniquely determined holomorphic
map $\tau:Y\mapsto Z$ such that $\tau \circ \Phi = \sigma$. The
map $\Phi$ blows down the zero section of $L'$. Note that the
set $A=\Phi(L'_{x_0}\cup X)= \Phi(L'_{x_0})$ is analytic in $Y$
(by Remmert's Proper Mapping Theorem) and the function
$\tilde{f}(\Phi(t)):=f(t)$ is holomorphic on $A$ (by property
(ii) of Remmert's reduction). Every analytic set in a complex
space is the support of a closed complex subspace (cf.
\cite{GR}, A.3.5), so we can apply Theorem V.4.4 in \cite{GR}
to conclude that the function $\tilde{f}$ is the restriction to
$A$ of a function $\tilde{F}$ that is holomorphic on the Stein
space $Y$. By the properties (ii) and (iii) above, there exists
a function $F$ holomorphic on $\Omega$ such that $\tilde{F}
\circ \Phi = F$. For $t\neq 0$ one can represent $F$ as $F\circ
\theta_i^{-1}(x,t)=\sum_{n=1}^{\infty}F_n^{(i)}(x)t^n$, with
$F_n^{(i)}$ holomorphic in $U_i$. We have $F\circ
\theta_k^{-1}(x,t)= F\circ \theta_i^{-1}(x,G_{ik}(x)t)$, which
gives $F_n^{(i)}(x)=(g_{ik}(x))^nF_n^{(k)}(x)$, i.e., $F_n$ are
cocycles corresponding to holomorphic sections of the tensor
product $L^n$ over $\Omega_r$. Considering the domain of
convergence of the representation for $F\circ \theta_k^{-1},
\quad k=1,...,l$, we get $\limsup_{n\to \infty}|F_n(x)|^{1/n}
\leq \exp rh(x), x \in X$. Let $\delta >0$. By Hartogs's lemma,
there exists an $n_\delta>1$ such that $(1/n)\log |F_n(x)|\leq
r\cdot h(x)+\delta, \quad x\in K, n\geq n_\delta$. For the
estimate from below, note that $F_n(x_0)=\chi_r(\zeta_0)=
rh(x_0)$ for all $n$. Since $rh=r(v+\varphi)$ is continuous,
there exists an $n_0 \geq n_\delta$ and a neighborhood
$W_{x_0}$ of $x_0$ such that
$(1/n_0)\log|F_{n_0}(x)|>rh(x)-\delta, \quad x \in W_{x_0}$.
Compactness of $K$ and suitable relations between
$\varepsilon,\delta$ and $r$ give holomorphic sections
satisfying the conclusion of the lemma.
\end{proof}

\begin{proof} {(of Theorem \ref{thm:reggrowth})}:
We mimic the method of proof of Theorem 2.8i in \cite{Bl1}. Let $u
\in PSH(X,\omega), u\mid_K \leq Q$. By Theorem 7.1 in \cite{GZ1},
there is a sequence $u_k \in PSH(X,\omega) \cap
\mathcal{C}^{\infty}(X)$ such that $u_k \searrow u$. Let
$\varepsilon >0$. By Dini's theorem, there exists an integer $k_0$
such that $u(x) \leq u_k(x) \leq Q(x) +\varepsilon$ for all $x \in
K, k\geq k_0$. By adding a small multiple of a local K\"{a}hler
potential in some coordinate neighborhood, we can assume that
$dd^cu_k+\omega$ is strictly positive. By the Approximation Lemma,
$\exists s_j^{(k)} \in \Gamma(X,L^{N_j^{(k)}}), \quad j=1,...,m_k$
such that
$$u_k-3\varepsilon \leq
\sup_{j=1,...,m_k}(\log|\exp(-2N_j^{(k)}\varepsilon
s_j^{(k)}|)/(N_j^{(k)})\leq (1/n)\log\Phi_n(x),$$ where $n =
\max_j N_j^{(k)}, j=1,...,m_k$. Hence $u-4\varepsilon \leq \log
\Phi$. The reverse inequality is obvious, since
$(1/N)\log\|s\|_{N\varphi}$ defines a positive singular metric
on $L$.
\end{proof}

Under the assumptions of Theorem \ref{thm:reggrowth} we also have the following:\\
\begin{prop}\label{prop:polygrowth}: Let
$\Psi(x)=\lim_{n\to\infty}\psi_n(x)=\sup_{n\geq 1}\psi_n(x),$
with $ \psi_n(x)= \sup\{\|s\|_{n\varphi}(x), s \in
\Gamma(X,L^n), \sup^\circ_K(\exp(-nQ)\|s\|_{n\varphi})\leq 1\}$
and $\sup^\circ_K(f) :=\inf\{\sup_{K\setminus P}(f): P\subset
K, P \mbox{ is }PSH(X,\omega)-\mbox{polar}\}$. Then
$$V_{K,\omega,Q}^* =(\log \Psi_K)^*.$$
\end{prop}
The proof proceeds exactly like that of \cite{Bl1}, Theorem
2.8(ii), provided we have the domination principle on a compact
K\"{a}hler manifold of dimension $N$ (cf. \cite{Kl}, cor. 3.7.5
and prop.5.5.1 \cite{BT2}, cor.4.5, \cite{Ta}, for versions on
open subsets of $\mathbb{C}^N$). In our proof we will assume that one of functions is in $L^{\infty}(X)$, since this is the case we need. A more general version was recently proved independently as Proposition 2.7 in \cite{BB}. Proofs of the domination principle rely on the comparison
principle, which was established in \cite{GZ2}(cf. also
\cite{Ko},\cite{HKH}) for the  class of functions
$\mathcal{E}(X,\omega)$ defined therein, which contains $L^{\infty}(X)$.  Recall the following result, which  allows us to 
apply the comparison and domination principles in the
weighted theory
\begin{prop}:  If $K$ is not $PSH(X,\omega)$-polar and $Q$ is continuous, then
$V_{K,\omega,Q}^*\in PSH(X,\omega)\cap L^{\infty}(X)$. In
particular, the complex Monge-Ampere operator
$(\omega_{V_{K,Q}})^n$ is well defined and satisfies
$(\omega_{V_{K,Q}})^N = 0$ in $X\setminus \overline{K}$.
\end{prop}
\begin{proof} The proof proceeds as that of \cite{GZ1}, Theorem 4.2.2,
and uses Proposition \ref{prop:nonpolar}.
\end{proof}

Now we may state and prove the required domination principle.

\begin{thm}(Domination Principle): Let $u,v \in
PSH(X,\omega)$ with $v \in L^{\infty}(X)$ be such that
$$\int_{\{u<v\}}(\omega+dd^cu)^N=0.$$ Then $u\geq v$ in X.
\end{thm}
\begin{proof} The following argument was communicated to us by Ahmed Zeriahi as a replacement for an earlier incorrect proof. It is enough to prove that $u \geq v$ on a set of full $\omega$-volume in $X$. We can assume that v is negative everywhere on $X$. Then for all $s,t>0$, $\{u-v \leq -s-t\} \subset \{u-v\leq -s-tv\}$, which for small $t$ is still a subset of $\{u-v<0\}$. Then, by Lemma 2.2 in \cite{EGZ},
$$ 0=\int_{\{u-v<-s-tv\}}(\omega + dd^cu)^N \geq t^N\mbox{Cap}\{u-v \leq -s-t\},$$ where Cap is the Monge-Ampere capacity defined in \cite{GZ1} (Definition 2.4). Proposition 2.5(1) in \cite{GZ1} implies that $\mbox{Vol}\{u-v \leq -s-t\}$ for $s,t>0$, $t$ small, hence $\mbox{Vol}\{u-v<0\}=0$.  
\end{proof}

Finally, we are interested in how weighted pluricomplex Green
functions change under a holomorphic map $f:X\to X$, where $X$
is a compact K\"{a}hler manifold (not necessarily projective
algebraic) with a closed real $(1,1)$-current $\omega$ on $X$ with
continuous local potentials (not necessarily a K\"{a}hler form).
Proposition 4.4.5 in \cite{GZ1} states that if $f:X \to X$ is
holomorphic, and $K\subset X$ is a Borel set, then
$V_{f(K),\omega}\circ f \leq V_{K,f^*\omega}$. The proof applies
also to the weighted pluricomplex Green
function and gives the following:\\
\begin{prop} Let $X,\omega, K$ be as above and let $Q$ be a mild
function on $X$. Then $V_{f(K),\omega,Q}\circ f \leq
V_{K,f^*\omega,Q\circ f}$ in $X$.
\end{prop}
Below, we establish a relation between the pullback of
$V_{K,\omega,Q}$ by a surjective holomorphic map $f:X\mapsto X$
and $V_{f^{-1}(K),\omega,\tilde{Q}}$ with an appropriate function
$\tilde{Q}$. For a  function $u:X \mapsto
\mathbb{R}\cup\{-\infty\}$ let us define $f_*u(x)=\sup \{u(y):y
\in f^{-1}(x)\}.$ This is a well defined function, since
$f^{-1}(x)$ is compact. Also, let $f^*u=u\circ f$. The following
theorem generalizes Theorem 2.12 in \cite{Bra}
and Theorem 1 in \cite{St1} (it yields both as special cases):\\

\begin{thm}\label{thm:pullback}: Assume that there exist $\alpha$
and $\beta$, $1<\alpha\leq \beta$,
such that $$\alpha
f_*(PSH(X,\omega))\subset PSH(X,\omega)$$ and
$$f^*(PSH(X,\omega))\subset \beta\cdot PSH(X,\omega).$$ Then for
every Borel set $K\subset X$ and every mild function $Q$ on $X$,
\[
\alpha V_{f^{-1}(K),\omega,f^*Q/\alpha}(x) \leq
V_{K,\omega,Q}\circ f(x) \leq \beta
V_{f^{-1}(K),\omega,f^*Q/\beta}.
\]
\end{thm}
\begin{proof} Let $u \in
PSH(X,\omega)$ be such that $\alpha u \leq f^*Q$ on $f^{-1}(K)$. Then
$v = \alpha f_*u$ is in $PSH(X,\omega)$ and satisfies $v \leq Q$ on
$X$. Moreover, $\alpha u(x) \leq v(f(x)) \leq
V_{K,\omega,Q}(f(x))$, which gives the first inequality. For the
second one, if $u \in PSH(X,\omega)$ satisfies $u \leq Q$ on $K$,
then by assumption $(1/\beta)f^*u$ is in $PSH(X,\omega)$ and
$(1/\beta)f^*u\leq (1/\beta)f^*Q$ on $f^{-1}(K)$, which gives the
conclusion.
\end{proof}

On $X = \mathbb{CP}^N$, the assumptions of Theorem
\ref{thm:pullback} are equivalent to assumptions about growth of
$f$ made in Theorem 2.12 in \cite{Bra} or its unweighted
counterpart, Theorem 5.3.1 in \cite{Kl}. Details may be found
in Theorem 1 in \cite{St1}) and its proof.  The main theorem in
\cite{St2} has conditions equivalent to the assumption
$\alpha f_*PSH(X,\omega)\subset PSH(X,\omega)$ when $X
\hookrightarrow \mathbb{CP}^N$ is a projective algebraic manifold
and $\omega$ is the pullback of the Fubini-Study form by the
embedding $\hookrightarrow$. One of the conditions is that $f$ has
an attracting
divisor in $X$, so in fact the assumption is quite strong.\\

\textbf{Acknowledgments:} The authors thank Thomas Bloom,
Finnur L\'{a}russon, Ragnar Sigurdsson and Ahmed Zeriahi for helpful comments onprevious drafts of this paper. In addition, the second named
author
 thanks the DFG Research Center Matheon, in particular Professor
Volker Mehrmann, for making it possible for her to stay and work
in Berlin in the academic year 2007/08.

\end{document}